\RequirePackage{fix-cm}

\documentclass[smallextended]{svjour3}       
\smartqed  %
\usepackage{amssymb}
\usepackage{graphicx}
\usepackage[all]{xy}
\usepackage[centertags]{amsmath}
\usepackage{latexsym}
\usepackage{amsfonts}
\usepackage{amssymb}

\begin{document}

\title{Computing the Theoretical Cost of the Optimal Ate Pairing on Elliptic Curves with Embedding Degrees 54 and 48 at the $256$-bits security level}
\titlerunning{Pairings with Elliptic Net}        
\author{Narcisse Bang Mbiang  \and   Emmanuel Fouotsa
    \and Celestin Lele}
\institute{ Narcisse Bang Mbiang\at
    Department of Mathematics and Computer Science,\\ Faculty of Sciences, University of Dschang (Cameroon)\\
    \email{$bang\_narcisse@yahoo.fr$}
    \and
     Emmanuel Fouotsa \at
    Department of Mathematics,~Higher Teacher Training
    College,\\ The University of Bamenda P.O BOX 39
    Bambili,~Cameroon.\\
    \email{emmanuelfouotsa@yahoo.fr}\\
   \and
    Celestin Lele \at
    Department of Mathematics and Computer Science,\\ Faculty of Sciences, University of Dschang (Cameroon)\\
    \email{$celestinlele@yahoo.fr$}
}

\date{Received:  / Accepted: }
\maketitle

\begin{abstract} Efficient computation of pairings with Miller
algorithm has recently received a great attention due to many
applications in cryptography. In this work,~we give formulae for the
optimal Ate pairing in terms of elliptic nets associated to twisted
Barreto-Naehrig (BN) curve, Barreto-Lynn-Scott(BLS) curves and
Kachisa-Schaefer-Scott(KSS) curves considered at the 128, 192 and
256-bit security levels. We show how to parallelize the computation
of these pairings when the elliptic net approach is used and we
obtain more efficient theoretical results with 8 processors compared
to the Miller loop approach for each corresponding case.
\end{abstract}

\section{Introduction}\label{intro}
$\hspace*{1cm}$Many new protocols such as the Identity-Based
Encryption \cite{BonFra01}, the tripartite Diffie-Hellman key
exchange \cite{Joux00} and short signatures \cite{BonLynSha04} are
based on pairings and so the efficiency of pairing computation has
become a field of active research.~The classical method for
computing pairings is the Miller's algorithm \cite{Miller04}. In
2007,~Katherine Stange \cite{Stange07} introduced a new algorithm to
evaluate pairings.~The algorithm is based on elliptic nets which is
the generalization of elliptic divisibility sequences to a higher
rank. Both methods compute pairings using $\mathcal{O}(\ell
og_{2}(r))$ field operations over an $r-$torsion subgroup. Based on
present results, the Miller algorithm remains the fastest method for
computing pairings. However, one can observes that formulae for the
Doubling and Addition steps in the  elliptic nets algorithm are
suitable for parallel calculations. This work aims at parallelizing
elliptic net algorithms for the optimal Ate pairing on
Barreto-Naehrig(BN) curve, Barreto-Lynn-Scott(BLS) curves with
embedding degree 12, 24 and 48 (\cite{BaLynScot03}) and
Kachisa-Schaefer-Scott(KSS) curves with embedding degree 16
(\cite{KaShEdSco08}), in other to compare their computational
costs to the Miller loop ones.\\

\textbf{Our contribution.}
\begin{enumerate}
    \item We give explicit
    formulae for computing the optimal Ate pairing on the above
    mentioned curves in terms of elliptic nets associated the
    twisted
    curves.
    \item We provide algorithms to parallelize the computation of
    the obtained elliptic net formulae.
    \item We give the computational
    cost for the part without the final exponentiation of optimal Ate
    pairing on each chosen curve in this work, and we compare these
    costs to those for the corresponding Miller loops.
\end{enumerate}
$\hspace*{1cm}$This paper is organized as follows. Section
\ref{sectionreview} summarizes notions on elliptic nets and elliptic
curves, Section \ref{sectionoptimalate} gives calculations of the
optimal Ate pairing on twisted elliptic curve in terms of elliptic
nets, Section \ref{sectionparallel} provides parallel algorithms for
$4$ and $8$ processors, the computational costs for the step without
the final exponentiation of optimal Ate pairings, and their
comparisons to those for the corresponding Miller loop, Section
\ref{conclusion} concludes the work.

%

\section{Elliptic nets and elliptic curves.}\label{sectionreview} $\hspace*{1cm}$In this section,~we define an elliptic net,~we
state the theorem which gives the bijection among elliptic nets and
elliptic curves and we express pairings in terms of elliptic nets.
\begin{definition}\label{def1}{(\cite{Stange07})}~Let~$\mathcal{A}$ be a finite generated free
    abelian group and~$\mathcal{R}$ an integral domain.~An elliptic net
    is a map~$W~:~\mathcal{A}\rightarrow\mathcal{R}$ which satisfies the
    recurrence relation
    $$W(p+q+s)W(p-q)W(r+s)W(r)+$$
    $$W(q+r+s)W(q-r)W(p+s)W(p)+$$
    \begin{equation}\label{eq3}
    W(r+p+s)W(r-p)W(q+s)W(q)=0
    \end{equation}
    for~$p,~q,~r,~s\in\mathcal{A}$
\end{definition}
Let~~$\{e_{i}\}_{i\in\{1,\cdots,n\}}$,~be the natural basis
of~$\mathcal{A}$. $W$ is normalized if~$W(e_{i})=1$ for all~$i$
and~$W(e_{i}+e_{j})=1$ for all~$1\leq i<j\leq n$.\\
$W$ is non-degenerate~if~$W(e_{i})\neq0$,~~$W(e_{i}+e_{j})\neq
0$,~$W(e_{i}-e_{j})\neq0$,~ ~$W(2e_{i})\neq0$.

\begin{theorem}\label{theo2}\cite[Page 54]{Stange02}
    Let~$W:\mathbb{Z}^{2}\rightarrow\mathcal{A}$ be a normalized
    non-degenerate elliptic net.~Then there is a curve~$E$ given by
    $$E~:~Y^{2}+a_{1}XY+a_{3}Y=X^{3}+a_{2}X^{2}+a_{4}X+a_{6}$$
    where
    $$a_{1}=\frac{W(2,0)-W(0,2)}{W(2,1)-W(1,2)}$$
    $$a_{2}=2W(2,1)-W(1,2)$$
    $$a_{3}=W(2,0),~~~~a_{4}=(W(2,1)-W(1,2))W(2,1),~~~~a_{6}=0$$
    with~$\psi(P,Q,\mathcal{C}_{ns})=W$,~where~$P=(0,0)$
    and~$Q=(W(1,2)-W(2,1),0)$ are non-singular points
    and~$\mathcal{C}_{ns}$ the non-singular part of the curve~$E$.
\end{theorem}

\begin{theorem}\cite{Stange07}
    Let $E$ be an elliptic curve defined over $\mathbb{C}$, and let
    $\Gamma$ be its corresponding lattice. Let $P_{1}$ and $P_{2}$ be
    two points in $E(\mathbb{C})$ such that $P_{1},~P_{2}\neq
    \mathcal{O}$ and $P_{1}\pm P_{2}\neq \mathcal{O}$. Let
    $z=(z_{1},z_{2})\in \mathbb{C}^{2}$ be such that $z_{1}$ and $z_{2}$
    correspond to $P_{1}$ and $P_{2}$ respectively, under the
    isomorphism $\mathbb{C}/\Gamma\cong E(\mathbb{C})$. For
    $v(v_{1},v_{2})\in \mathbb{Z}^{2}$, define a function
    $\psi_{(v_{1},v_{2})}$ on $\mathbb{C}^{2}$ in variables
    $z=(z_{1},z_{2})$ as follows:
    $$\psi_{(v_{1},v_{2})}(z_{1},z_{2};~\Gamma)=\frac{\sigma(v_{1}z_{1}+v_{2}z_{2};~\Gamma)}
    {\sigma(z_{1};~\Gamma)^{v^{2}_{1}-v_{1}v_{2}}\sigma(z_{1}+z_{2};~\Gamma)^{v_{1}v_{2}}\sigma(z_{2};~\Gamma)^{v^{2}_{2}-v_{1}v_{2}}},$$
    we have $$\psi_{v}(P_{1},P_{2};~E(\mathbb{C}))
    =\psi_{(v_{1},v_{2})}(z_{1},z_{2};~\Gamma)$$ and
\end{theorem}
\[\begin{array}{ccccc} W: & \mathbb{Z}^{2}
& \to & \mathbb{C} & \\
& (v_{1},v_{2}) &\mapsto & \psi_{(v_{1},v_{2})}(z_{1},z_{2};~\Gamma),\\
\end{array}\]\\
is an elliptic net associated to the curve $E$ and the points
$P_{1}$ and $P_{2}$.\\
Stange proved this result to any field $\mathbb{K}$ \cite{Stange07}.
In this work,~we consider an elliptic curve~$E$ in the reduced
form~$y^{2}=x^{3}+Ax+B$,~$P(x_{1},y_{1}),~Q=(x_{2},y_{2})\in
E(\mathbb{F}_{p^{k}})$. Initial values of the elliptic nets $W(i,0)$
and $W(i,1)$ associated to $E, P, Q$ are:
\begin{eqnarray}\label{equan 1}
W(1,0)&=& 1,
\end{eqnarray}
\begin{eqnarray}\label{equan 2}
W(2,0) &=& 2y_{1},
\end{eqnarray}
\begin{eqnarray}\label{equan 3}
W(3,0)&=& 3x^{3}_{1}+6Ax^{2}_{1}+12Bx_{1}-A^{2},
\end{eqnarray}
\begin{eqnarray}\label{equan 4}
W(4,0)&=&
4y_{1}(x^{6}_{1}+5Ax^{4}_{1}+20Bx^{3}_{1}-5A^{2}x^{2}_{1}-4ABx_{1}-8B^{2}-A^{3}),
\end{eqnarray}
\begin{eqnarray}\label{equan 5}
W(0,1) &=&W(1,1)=1,
\end{eqnarray}
\begin{eqnarray}\label{equan 6}
W(2,1)&=& 2x_{1}+x_{2}-(\frac{y_{2}-y_{1}}{x_{2}-x_{1}})^{2}
\end{eqnarray}
\begin{eqnarray}\label{equan 7}
W(-1,1) &=& x_{1}-x_{2} ,
\end{eqnarray}
\begin{eqnarray}\label{equan 8}
W(2,-1)&=& (y_{1}+y_{2})^{2}-(2x_{1}+x_{2})(x_{1}-x_{2})^{2}.
\end{eqnarray}
There is a bijection between the set of elliptic curves and two
points $P$ and $Q$ with $P$, $Q$, $P+Q$ and $P-Q\neq \mathcal{O}$,
and the set of elliptic nets $W(n,m)$ associated to $E, P, Q$, such
that $W(1,0)=W(0,1)=W(1,1)=1$ and $W(1,-1)\neq 0$. This works for
higher ranks. Since pairings are classically defined on elliptic
curves,~one can then express
pairings in terms of elliptic nets.\\
\subsection{Tate Pairing in terms of Elliptic Nets.}
$\hspace*{1cm}$In this subsection, we are summarizing Stange
construction of the Tate pairing via elliptic net and also the
double-and-add algorithm for its computation \cite{Stange07}.\\
Let~$E$ be an elliptic curve defined over a finite
field~$\mathbb{F}_{p}$, with $p$ a prime number greater than 3. Let
$r$ be a large prime dividing the order of the elliptic curve group
such that $gcd(p,r)=1$. Let~$k$ be the smallest positive integer and
also called an embedding degree of the curve with respect to $r$.
Let $P\in E(\mathbb{F}_{p})[r]$ and~$Q\in
E(\mathbb{F}_{p^{k}})[r]$.~The Tate paring is
\begin{equation}\label{eq4}
T_{r,P}=\frac{W_{P,Q}(r+1,1)W_{P,Q}(1,0)}
{W_{P,Q}(r+1,0)W_{P,Q}(1,1)}
\end{equation}
where~$W_{P,Q}$ is the elliptic net associated to~$E,~P$ and~$Q$.\\
One can show from Definition 1 that~$W_{P,Q}(0,0)=0$
and~$W_{P,Q}(-a,-b)=-W_{P,Q}(a,b)$. Since elliptic nets which are in
bijection with  elliptic curves
satisfy~$W_{P,Q}(1,0)=W_{P,Q}(1,1)=W_{P,Q}(1,1)=1$,~the Tate pairing
is then
\begin{equation}\label{eq5}
T_{r,P}=\frac{W_{P,Q}(r+1,1)} {W_{P,Q}(r+1,0)}
\end{equation}
Stange provided a double-and-add algorithm for
computing~$W_{P,Q}(r+1,0)$ and~$W_{P,Q}(r+1,1)$ in
~$log_{2}(r+1)-1$ steps.\\
The method consists on defining an initial block~$V$ (Table 1) with
a first vector of 8 consecutive terms of the sequence~$W(i,0)$
centered on~$W(k,0)$ and~$W(k+1,0)$ and a second vector of 3
consecutive terms of the sequence~$W(i,1)$ centered on the
term~$W(k,1)$.~Two
functions \textbf{Double(V)} and \textbf{DoubleAdd(V)} are provided and described as follows:\\
$\textbf{Double(V):}$~given a block~$V$ centered at~$k$,~returns the
block centered at~$2k$.\\
$\textbf{DoubleAdd(V):}$~given a block~$V$ centered at~$k$ ,~returns
the block centred at~$2k+1$.\\\\
\begin{table}\centering
\begin{tabular}{|c|c|c|c|c|c|c|c|}
\multicolumn{1}{c}{}&\multicolumn{1}{c|}{}&(k-1,1)&(k,1)&(k+1,1)&\multicolumn{1}{c}{}&\multicolumn{1}{c}{} \\
\hline
(k-3,0)&(k-2,0)&(k-1,0)&(k,0)&(k+1,0)&(k+2,0)&(k+3,0)&(k+4,0)\\
\hline
\end{tabular}\caption{Block V centered at k}.
\end{table}

Based on Definition $\ref{def1}$, \textbf{Double(V)} and
\textbf{DoubleAdd(V)} can be obtained from the following
proposition:
\begin{proposition}\label{prop1}\cite{Stange07} Let $W$ be an elliptic net
    associated to an elliptic curve $E$ and 2 rational points. We then
    have the following relations: {\scriptsize
        \begin{flushleft}
            $W(2i-1,0)=W(i+1,0)W(i-1,0)^{3}-W(i-2,0)W(i,0)^{3},$
    \end{flushleft}}

    \begin{flushleft}
        $W(2i,0)=\frac{1}{W(2,0)}(W(i,0)W(i+2,0)W(i-1,0)^{2}-W(i,0)W(i-2,0)W(i+1,0)^{2}),$
    \end{flushleft}

    for~$i=k-1\cdots,~k+3$ and {\scriptsize
        \begin{flushleft}
            $W(2k-1,1)=\frac{1}{W(1,1)}(W(k+1,1)W(k-1,1)W(k-1,1)^{2}-W(k,0)W(k-2,0)W(k,1)^{2}),$
        \end{flushleft}
        \begin{flushleft}
            $W(2k,1)= W(k-1,1)W(k+1,1)W(k,0)^{2}-W(k-1,0)W(k+1,0)W(k,1)^{2},$
        \end{flushleft}
        \begin{flushleft}
            $W(2k+1,1)=\frac{1}{W(-1,1)}(W(k-1,1)W(k+1,1)W(k+1,0)^{2}-W(k,0)W(k+2,0)W(k,1)^{2}),$
        \end{flushleft}
        \begin{flushleft}
            $W(2k+2,1)=\frac{1}{W(-2,1)}(W(k-1,1)W(k+1,1)W(k+2,0)^{2}-W(k+1,0)W(k+3,0)W(k,1)^{2}).$
        \end{flushleft}
    }
\end{proposition}

\subsection{A Simplified Tate pairing in terms of Elliptic Net.}
$\hspace*{1cm}$N.~Ogura,~N.~Kanayama,~S.~Uchiyama,~and E.~Okamoto
\cite{NaoKaSka11} have expressed the reduced Tate pairing in term of
elliptic net.~They used a process of normalization for elliptic net
function to simplify the Tate pairing as it is given in the
following definition.
\begin{definition}{(\cite{NaoKaSka11})}~Let~$E$ be an elliptic curve defined over a
    finite field~$\mathbb{F}_{p}$,~$r$ a large prime number such
    that~$r$ divides~$\sharp E(\mathbb{F}_{p})$ and~$gcd(r,p)=1$.~$k$
    the embedding degree of the curve~$E$.~Let~$P\in
    E(\mathbb{F}_{p})[r]$ and~$Q\in E(\mathbb{F}_{p^{k}})[r]$,
    \begin{eqnarray}\label{eq99}
    T^{Red}_{r,P}(P,Q) &=& f_{r,P}(Q)^{\frac{p^{k}-1}{r}}=
    W_{P,Q}(r,1)^{\frac{p^{k}-1}{r}}
    \end{eqnarray}
    where~$W_{P,Q}$ is an elliptic net associated to~$E,~P~and~Q$.
\end{definition}
\section{Computation of the Optimal Ate Pairings on Twisted elliptic Curves in terms of Elliptic Nets} \label{sectionoptimalate}
Pairing-friendly curves are generally parameterized as $(p,r,t)$
where $p,r$ are given as polynomials in a variable $x$ and $t$ is
the trace of the Frobenius of the curve. A value of $x$ gives the
suitable primes $p$ and $r$ defining an elliptic curve with
cardinality $p+1-t$ divisible by $r$ at the corresponding security
level. This section describes optimal ate pairing on
pairing-friendly elliptic curves in terms of elliptic nets.
\subsection{The Optimal Ate Pairing}
$\hspace*{1cm}$The method for the construction of the optimal Ate
pairing is given in \cite{Ver10}. Let $\pi_{p}$ be the Frobenius map
on an elliptic curve defined by $\pi_{p}(x,y)=(x^{p},y^{p})$.  Let
$t$ be the trace of the Frobenius on $E(\mathbb{F}_{p})$ and
$T=t-1$. Let $P\in
\mathbb{G}_{1}=E(\mathbb{\overline{F}}_{p})[r]\cap Ker(\pi_{p}-[1])$
and $Q\in\mathbb{G}_{2}= E(\mathbb{\overline{F}}_{p})[r]\cap
Ker(\pi_{p}-[p])$, that means $Q$ satisfies $\pi_{p}(Q)=[p]Q$.\\
Let $\ell=mr$ be a multiple of $r$ such that $r\nmid m$ and write
$\ell=\sum_{i=0}^lc_ip^i=h(p)$, $(h(z)\in \mathbb{Z}[z])$. Recall
that $h_{R,S}$ is the Miller function \cite{Miller04}. For
$i=0,\cdots l$ set $s_i=\sum_{j=i}^lc_jp^j$; then the map
\begin{eqnarray}
\begin{array}{cccc}\label{op}
e_o:&\mathbb{G}_{2}\times \mathbb{G}_{1}& \rightarrow &\mu _{r}\\
&(Q,P)&\longmapsto & (\prod_{i=0}^lf_{c_i,Q}^{p^i}(P)\cdot\\
&&&\prod_{i=0}^{l-1}h_{[s_{i+1}]Q,[c_ip^i]Q}(P))^{\frac{p^k-1}{r}}
\end{array}
\end{eqnarray}
defines a bilinear pairing and non degenerate if \\$mkp^k \neq
((p^k-1)/r)\cdot \sum_{i=0}^lic_ip^{i-1} \mbox{ mod } r.$ The
coefficients $c_i: i=0,\cdots,l$ can be obtained from the short
vectors obtained from the lattice
\begin{eqnarray}\label{lattice}
L=\left(\begin{array}{ccccc}
r & 0 & 0 &\cdots &0 \\
-p & 1 & 0 & \cdots &0 \\
-p^2 & 0 & 1 & \cdots &0 \\
\cdots &  \cdots &  \cdots &  \cdots &  \cdots \\
-p^{\varphi(k)-1} & 0 & 0 & \cdots &1 \\
\end{array}
\right)
\end{eqnarray}

\begin{theorem}\label{theo1}\cite[Page 71]{Stange02}.
    Let $E$ be an elliptic curve defined over a field
    $\mathbb{F}_{p^{k}}$. Let $E_{\theta}$ be the twist curve of $E$
    defined over $\mathbb{F}_{p^{k/\delta}}$, where $\delta$ is the
    degree of twist, $\theta$ a generator of the basis of
    $\mathbb{F}_{p^{k}}$, seen as a $\mathbb{F}_{p^{k/\delta}}-$vector
    space. Let $\sigma_{\theta}:~E_{\theta}\rightarrow E$ be the
    twisting isomorphism such that
    $(x,y)\mapsto(x\theta^{2},y\theta^{3})$, $E:~y^{2}=x^{3}+b$ and
    $E_{\theta}:~y^{2}=x^{3}+b\theta^{-6}$.
    Let~$W_{\widetilde{Q},\widetilde{P}}$ be the elliptic net associated
    to a twist curve~$E_{\theta}$ and the
    points~$\widetilde{Q},~\widetilde{P}$, such that
    $\sigma_{\theta}(\widetilde{P})=P$ and
    $\sigma_{\theta}(\widetilde{Q})=Q$ respectively. We then have the
    following relations:
    \begin{equation}\label{eq1000}
    W_{Q,P}(n,0)=\theta^{1-n^{2}}W_{\widetilde{Q},\widetilde{P}}(n,0)
    \end{equation}
    \begin{equation}\label{eq100}
    W_{Q,P}(n,1)=\theta^{n-n^{2}}W_{\widetilde{Q},\widetilde{P}}(n,1)
    \end{equation}
\end{theorem}
One can prove Theorem $\ref{theo1}$ by induction.\\
We can now express the optimal Ate pairing in terms of elliptic nets
associated to the twist curves.\\
In the following,  $\widetilde{W}$ denotes
$W_{\widetilde{Q},\widetilde{P}}$ and $W$ denotes $W_{Q,P}$.
\subsection{Pairing on Twisted BN-Curves}
$\hspace*{1cm}$The family of BN-curves \cite{BareNa06} has embedding
degree $k=12$ and is given by the following parametrization:
\begin{eqnarray*}
    p&=&36x^4+36x^3+24x^2+6x+1, \\
    r&=& 36x^4+36x^3+18x^2+6x+1\\
    t&=& 6x^2+1
\end{eqnarray*}
The optimal Ate pairing for BN-curves is given in \cite{Ver10}  by:
\[\begin{array}{ccccc} e_{1}: &\mathbb{G}_{2}
\times\mathbb{G}_{1}
& \to & \mu_{r} & \\
& (Q,P) &\mapsto & (f_{6x+2,Q}\cdot \ell_{[6x+2]Q,[p]Q}\cdot \ell_{[6x+2+p]Q,[-p^2]Q}(P))^{\frac{p^{12}-1}{r}},\\
\end{array}\]\\
where $f_{n,Q}$ is the Miller function \cite{Miller04} and
$\ell_{R,S}$, the line passing through $R$ and $S$. The optimal Ate
pairing in terms of elliptic nets associated to twisted BN-curve is
already calculated in \cite{HiroOnu16}. The costs of the Double and
DoubleAdd steps for the parallelization of the Elliptic Net
Algorithm is given for 1, 4, 6, 8 and 10 processors. We are focus on
4 and 8 processors in this work. We little bit improve the costs of
the \textbf{Double} and \textbf{DoubleAdd} steps for 4 and 8
processors compared to theirs and we also noticed that one can
appreciate the importance of the parallelization of the Elliptic Net
Algorithm compared to the Miller loop from these numbers of
processors. We will give the cost of the part without the final
exponentiation of the optimal Ate pairing in
terms of elliptic nets which is not yet given.\\
The BN-curve $j-$invariant 0 has a twist of order
$\delta=6$.~Let~$W$ be the elliptic net associated to the BN-curve
$E$ and the points~$Q,P$,~and~$\widetilde{W}$ the elliptic net
associated to the twist curve as considered in \cite{HiroOnu16}.\\

Using~$(\ref{eq100})$ and the fact that the final exponentiation
eliminate~$\theta$, The optimal Ate pairing is given by:\\
\[\begin{array}{ccccc} e_{1}: &\mathbb{G}_{2}
\times\mathbb{G}_{1}
& \to & \mu_{r} & \\
& (Q,P) &\mapsto & (\widetilde{W}(6x+2,1)\cdot\mathbb{L}_{1}\cdot\mathbb{L}_{2})^{\frac{p^{12}-1}{r}},\\
\end{array}\]
where $\mathbb{L}_{1}=\ell_{[6x+2]\widetilde{Q},[p]\widetilde{Q}}$
and
$\mathbb{L}_{2}=\ell_{[6x+2+p]\widetilde{Q},[-p^{2}]\widetilde{Q}}$
are the line evaluations. The following theorem helps to compute the
coordinates in terms of elliptic nets for a multiple point.

\begin{theorem}\label{theomulelli}\cite[page 11]{NaoKaSka11} Let $n\in \mathbb{Z}$ and $S=(x,y)$, a point on an elliptic curve
$E$. The multiple point $[n]S$ in terms of elliptic nets is given by:\\
\begin{equation}\label{mul-ellp-net}
   [n]S=(x-\frac{\Psi_{n-1}\Psi_{n+1}}{\Psi^{2}_{n}}(x,y),\frac{\Psi^{2}_{n-1}\Psi_{n+2}-\Psi^{2}_{n+1}\Psi_{n-2}}{4y\Psi^{3}_{n}}(x,y)),
\end{equation}
where $\Psi_{n}(x,y)$ is the elliptic net of rank 1 associated to
$[n]S$
\end{theorem}

These lines are calculated as follow:\\

Coordinates of $[6x+2]\widetilde{Q}$ :\\\\
$x_{[6x+2]\widetilde{Q}}=\frac{S}{\widetilde{W}(6x+2,0)^{2}}$ where $S=x_{\widetilde{Q}}\widetilde{W}(6x+2,0)^2-\widetilde{W}(6x+1,0)\widetilde{W}(6x+3,0)$\\\\
$y_{[6x+2]\widetilde{Q}}=\frac{T}{\widetilde{W}(6x+2,0)^{3}}$ where
$T=\frac{\widetilde{W}(6x+1,0)^2\widetilde{W}(6x+4,0)-\widetilde{W}(6x+3,0)^{2}\widetilde{W}(6x,0)}{4y_{\widetilde{Q}}}$.\\\\
Coordinates of $[p]\widetilde{Q}$ :\\\\
$x_{[p]\widetilde{Q}}=\theta^{2(p-1)}x^{p}_{\widetilde{Q}}$,~~$y_{[p]\widetilde{Q}}=\theta^{3(p-1)}y^{p}_{\widetilde{Q}}$.\\\\

$\frac{y_{[6x+2]\widetilde{Q}}-y_{[p]\widetilde{Q}}}{x_{[6x+2]\widetilde{Q}}-x_{[p]\widetilde{Q}}}=\frac{T-\theta^{3(p-1)}\widetilde{W}(6x+2,0)^{3}
y^{p}_{\widetilde{Q}}}{S-\theta^{2(p-1)}\widetilde{W}(6x+2,0)^{2}
x^{p}_{\widetilde{Q}}}\times \frac{1}{\widetilde{W}(6x+2,0)}$~~,~~

$\frac{y_{\widetilde{P}}-y_{[6x+2]\widetilde{Q}}}{x_{\widetilde{P}}-x_{[6x+2]\widetilde{Q}}}=\frac{y_{\widetilde{P}}\widetilde{W}(6x+2,0)^{3}-T}
{x_{\widetilde{P}}\widetilde{W}(6x+2,0)^{2}-S}\times
\frac{1}{\widetilde{W}(6x+2,0)}$\\\\
Equation of the line evaluation $\mathbb{L}_{1}$ is then :\\\\
$\mathbb{L}_{1}:~(S-\theta^{2(p-1)}\widetilde{W}(6x+2,0)^{2}x^{p}_{\widetilde{Q}})(y_{\widetilde{P}}\widetilde{W}(6x+2,0)^{3}-T)-(T-\theta^{3(p-1)}\widetilde{W}(6x+2,0)^{3}y^{p}_{\widetilde{Q}})(x_{\widetilde{P}}\widetilde{W}(6x+2,0)^{2}-S)$\\\\

Coordinates of $[6x+2+p]\widetilde{Q}$ :\\\\
$[6x+2+p]\widetilde{Q}=[6x+2]\widetilde{Q}
+[p]\widetilde{Q}=(\frac{S}{\widetilde{W}(6x+2,0)^{2}}~,~\frac{T}{\widetilde{W}(6x+2,0)^{3}}
)+(\theta^{2(p-1)}x^{p}_{\widetilde{Q}}~,~\theta^{3(p-1)}y^{p}_{\widetilde{Q}})$\\\\
The calculation with good reduction using the equation of the
twisted curve gives:\\\\
$x_{[6x+2+p]\widetilde{Q}}=\frac{U}{Z^{2}}$
where $U=2b\theta^{-6}\widetilde{W}(6x+2,0)^{6}+S\widetilde{S}(S+\widetilde{S})-2T\widetilde{T}$.\\\\
$\widetilde{S}=\theta^{-2(p-1)}\widetilde{W}(6x+2,0)^{2}x^{p}_{\widetilde{Q}}$,~~~
$\widetilde{T}=\theta^{-3(p-1)}\widetilde{W}(6x+2,0)^{3}y^{p}_{\widetilde{Q}}$,~~~
$Z=(S-\widetilde{S})\widetilde{W}(6x+2,0)$.\\\\
$y_{[6x+2+p]\widetilde{Q}}=\frac{V}{Z^{3}}$ where
$V=(T-\widetilde{T})(T\widetilde{T}-3\theta^{-6}\widetilde{W}(6x+2,0)^{6}b)+3S\widetilde{S}(S\widetilde{T}-T\widetilde{S})$.\\\\
The terms $\widetilde{S},~\widetilde{T}$ and $Z$ are as defined
above.\\
Coordinates of $[-p^2]\widetilde{Q}$ :\\\\
$x_{[-p^2]\widetilde{Q}}=\theta^{2(p^2-1)}x^{p^2}_{\widetilde{Q}}=\theta^{2(p^2-1)}x_{\widetilde{Q}}$~~,~~$y_{[-p^2]\widetilde{Q}}=-\theta^{3(p^2-1)}y^{p^2}_{\widetilde{Q}}=-\theta^{3(p^2-1)}y_{\widetilde{Q}}$.\\\\
Since the coordinates of $\widetilde{Q}$ are in
$\mathbb{F}_{p^{2}}$.\\
Equation of the line evaluation $\mathbb{L}_{2}$ with good simplification is then :\\\\
$\mathbb{L}_{2}:~(y_{\widetilde{P}}Z^3-V)(U-
\theta^{2(p^2-1)}x_{\widetilde{Q}}Z^2)-(x_{\widetilde{P}}Z^2-U)(V+\theta^{3(p^2-1)}y_{\widetilde{Q}}Z^{3})$\\\\

%
%

\subsection{Pairing on Twisted BLS-Curves of Embedding degrees $12$, $24$ and $48$}
The $BLS12$ family of elliptic curves \cite{BaLynScot03} are
parameterized by:
\begin{eqnarray*}
    p&=&(x-1)^{2}(x^{4}-x^{2}+1)/3+x, \\
    r&=& x^{4}-x^{2}+1\\
    t&=& x+1
\end{eqnarray*}
The $BLS24$ family of elliptic curves \cite{BaLynScot03} are
parameterized by:
\begin{eqnarray*}
    p&=&(x-1)^{2}(x^{8}-x^{4}+1)/3+x, \\
    r&=& x^{8}-x^{4}+1\\
    t&=& x+1
\end{eqnarray*}
The $BLS48$ family of elliptic curves \cite{BaLynScot03} are
parameterized by:
\begin{eqnarray*}
    p&=&(x-1)^{2}(x^{16}-x^{8}+1)/3+x, \\
    r&=& x^{16}-x^{8}+1\\
    t&=& x+1
\end{eqnarray*}
The optimal Ate pairing for $BLS12$, $BLS24$ and $BLS48$ curves is
given in \cite{HesSmaVer06} by:
\[\begin{array}{ccccc} e_{2}: &\mathbb{G}_{2}
\times\mathbb{G}_{1}
& \to & \mu_{r} & \\
& (Q,P) &\mapsto & f_{x,Q}(P)^{\frac{p^{k}-1}{r}},\\
\end{array}\]\\
where $f_{x,Q}$ is the Miller function \cite{Miller04} and $k=12$,~$k=24$ and $k=48$ respectively.\\
$\hspace*{1cm}$The BLS12 BLS24 and BLS48 curves with $j-$invariant 0
have a twists of order $\delta=6$. From the isomorphism described in
theorem $\ref{theo1}$, let~$W$ be the elliptic net associated to the
BLS12, BLS24 or BLS48 curve and the
points~$Q,P$,~and~$\widetilde{W}$ the elliptic net associated to the
twist~$E_{\theta}$ and the points~$\widetilde{Q},~\widetilde{P}$,~
where~$\widetilde{Q},~\widetilde{P}$ correspond to~$Q,~P$ via
$\sigma_{\theta}$ respectively.~Using~$(\ref{eq100})$ and the fact
that the final exponentiation eliminates~$\theta$,~we have
\begin{equation*}
e_{2}=f_{x,Q}(P)^{\frac{p^{k}-1}{r}}=\widetilde{W}(x,1)^{\frac{p^{k}-1}{r}},
\end{equation*}
where $k=12,~24$ or $48$ respectively.
\subsection{Pairing on Twisted KSS Curves of Embedding Degree $16$}
$\hspace*{1cm}$The $KSS16$ family of elliptic curves
\cite{KaShEdSco08} are parameterized by:
\begin{eqnarray*}
    p&=&\frac{1}{980}(x^{10}+2x^9+5x^8+48x^6+152x^5+240x^4+625x^2+2398x+3125), \\
    r&=& \frac{1}{61250}(x^{8}+48x^{4}+625)\\
    t&=&\frac{1}{35}(2x^5+41x+35)
\end{eqnarray*}
The optimal Ate pairing for KSS16 curves is given in \cite{Ver10} by:\\
\[\begin{array}{ccccc} e_{3}: &\mathbb{G}_{2}
\times\mathbb{G}_{1}
& \to & \mu_{r} & \\
& (Q,P) &\mapsto & ((f_{x,Q}(P)\cdot
\ell_{[x]Q,[p]Q}(P))^{p^{3}}\cdot
\ell_{Q,Q}(P))^{\frac{p^{16}-1}{r}},\\
\end{array}\]\\\\
where $f_{x,Q}$ is the Miller function \cite{Miller04}, and where
$\ell_{[x]Q,[p]Q}$ is the line trough $[x]Q$ and $[p]Q$, and
$\ell_{Q,Q}(P)$, the tangent line trough $Q$.\\
The $j-$invariant is 1728 for KSS16-curve and then has a twist of
order 4.~Let~$W$ be the elliptic net associated to the KSS16-curve
$E:~y^{2}=x^{3}+ax$ and the points~$Q,P$,~and~$\widetilde{W}$ the
elliptic net associated to the quartic
twist~$E_{\eta}:~y^{2}=x^{3}+\eta^{-4}ax$ where
$\eta\in\mathbb{F}^{\star}_{p^{16}}$,~$(1,\eta,\eta^{2},
\eta^{3})$,~a basis of the~$\mathbb{F}_{p^{4}}$-vector
space~$\mathbb{F}_{p^{16}}$, and the
points~$\widetilde{Q},~\widetilde{P}$,~
where~$\widetilde{Q},~\widetilde{P}$ correspond to~$Q,~P$
respectively, via the isomorphism
~$\sigma_{\eta}:~E_{\eta}(\mathbb{F}_{p^{4}})\rightarrow
E(\mathbb{F}_{p^{16}}),~(x,y)\mapsto(\eta^{2}x,\eta^{3}y)$.\\
Using~$(\ref{eq100})$ and the fact that the final exponentiation
eliminates~$\theta$,~we have
\begin{equation*}
f_{x,Q}(P)^{\frac{p^{16}-1}{r}}=\widetilde{W}(x,1)^{\frac{p^{16}-1}{r}}
\end{equation*}
Let's compute $\ell_{1}$ and $\ell_{1}$ the values $\ell_{[x]\widetilde{Q},[p]\widetilde{Q}}(\widetilde{P})$ and $\ell_{2}= \ell_{\widetilde{Q},\widetilde{Q}}(\widetilde{P})$ respectively:\\
Based on Theorem $\ref{theomulelli}$, we then have:\\\\
Coordinates of $[x]\widetilde{Q}$ :\\\\
$x_{[x]\widetilde{Q}}=\frac{A}{\widetilde{W}(x,0)^{2}}$ where
$A=x_{\widetilde{Q}}\widetilde{W}(x,0)^2-\widetilde{W}(x-1,0)\widetilde{W}(x+1,0)$.\\\\
$y_{[x]\widetilde{Q}}=\frac{B}{\widetilde{W}(x,0)^{3}}$ where
$B=\frac{\widetilde{W}(x-1,0)^2\widetilde{W}(x+2,0)-\widetilde{W}(x+1,0)^{2}\widetilde{W}(x-2,0)}{4y_{\widetilde{Q}}}$.\\\\
Coordinates of $[p]\widetilde{Q}$ :\\\\
$x_{[p]Q}=\eta^{2(p-1)}x^{p}_{\widetilde{Q}}$,~~$y_{[p]Q}=\eta^{3(p-1)}y^{p}_{\widetilde{Q}}$.\\\\
Equations of the line evaluations $\ell_{1}$ and $\ell_{2}$:\\\\
$\ell_{1}:
(A-\eta^{2(p-1)}\widetilde{W}(x,0)^{2}x^{p}_{\widetilde{Q}})(\widetilde{W}(x,0)^{3}y_{P}-B)-
(B-\eta^{3(p-1)}\widetilde{W}(x,0)^{3}y^{p}_{\widetilde{Q}})(\widetilde{W}(x,0)^{2}x_{P}-A)$\\\\
$\ell_{2} :
(3\eta^{4}x^{2}_{\widetilde{Q}}+a)x_{\widetilde{P}}-2\eta^{4}y_{\widetilde{Q}}y_{\widetilde{P}}+2\eta^{4}y^{2}_{\widetilde{Q}}
-3\eta^{4}x^{3}_{\widetilde{Q}}-ax_{\widetilde{Q}}$.\\\\
We then have
\[\begin{array}{ccccc} e_{3}: &\mathbb{G}_{2}
\times\mathbb{G}_{1}
& \to & \mu_{r} & \\
& (Q,P) &\mapsto & ((\widetilde{W}(x,1)\cdot\ell_{1})^{p^{3}}\cdot\ell_{2})^{\frac{p^{16}-1}{r}},\\
\end{array}\]

$\hspace*{1cm}$ In Table 2, we give some good parameters collected
from \cite{RazDuq18} for the chosen BN, BLS12, KSS16, and BLS24
curves and we also provide a good parameter for BLS48 curve in order
to give the computational cost of the Miller
loop for each corresponding optimal Ate pairing.\\
\begin{table}
    \centering\scriptsize
    \begin{tabular}{|c|c|c|c|c|c|c|}
        \hline
        Sec. level & Curve & Curve Eq. & k & Parameter $x$ & $\lceil log_{2}r\rceil$  & $\lceil log_{2}p\rceil$\\
        \hline
        & BN\cite{RazDuq18} & $y^2=x^3-4$ &12 &$x=2^{114}+2^{101}-2^{14}-1$& 280 & 280 \\
        128-bit & BLS12\cite{RazDuq18} &$y^2=x^3+4$ &12& $x=-2^{77}+2^{50}+2^{33}$& 273& 616 \\
        & KSS16\cite{RazDuq18} & $y^2=x^3+x$ &16 &$x=2^{35}-2^{32}-2^{18}+2^8+1$& 281 & 340 \\
        \hline\hline

        192-bit & BLS24\cite{RazDuq18} & $y^2=x^3-2$ &24& $x=-2^{56}-2^{43}+2^{9}-2^6$ & 427 & 558 \\
        \hline\hline

        & BLS24\cite{RazDuq18}& $y^2=x^3-2$ &24 &$x=-2^{103}-2^{101}+2^{68}+2^{50}$ & 581 & 1028  \\
      256-bit   & BLS48& $y^2=x^3+11$ & 48&$x=2^{32}-2^{18}-2^{10}-2^4$ & 512 & 575  \\

        \hline
    \end{tabular}\caption{Selected Parameters for our chosen curves.}
\end{table}

\subsection{Computational Costs of the Miller Loop for the Studied curves}
 In this section, we give the computational costs for the
Miller loop \cite{Miller04} using the selected parameters in Table
2. Computational costs at $128$-bit security level are given in
\cite{RazDuq18}. That's, BN-curve ($12068M$), $BLS12$-curve
($7708M$) and $KSS16$-curve ($7534M$) . The remaining computational
costs for $192$-bit security level and $256$-security level
will be provided in this work.\\

The most efficient formulae for the doubling steps (doubling of
points and line evaluations) and the addition steps (addition of
points and line evaluations) for curves with sextic twists are given
in \cite{CosLanNae10}. The doubling step costs $53M$ and the
addition step is $76M$.



The best choice for the parameters $x$ at $192$-bit and $256$-bit
security levels for $BLS24$-curves given in \cite{RazDuq18} are
$x=-2^{56}-2^{43}+2^{9}-2^6$ and $x=-2^{103}-2^{101}+2^{68}+2^{50}$
respectively. The doubling step costs $68M$ and the addition step is
$110M$. For the parameter $x=-2^{56}-2^{43}+2^{9}-2^6$, the
computational cost for the Miller loop is 56 doubling steps, 3
addition steps, 55 squarings and 58 multiplications in
$\mathbb{F}_{p^{18}}$. That is, $56(68M)+55S_{24}+3(110M)+58M_{24}$.
That's $19474M$. For the parameter
$x=-2^{103}-2^{101}+2^{68}+2^{50}$, the computational cost for the
Miller loop is 103 doubling steps, 3 addition steps, 102 squarings
and 105 multiplications in $\mathbb{F}_{p^{18}}$. That is,
$103(68M)+102(108M)+3(110M)+105(162)=35360M$. In the same
consideration, the computational costs of the Miller loop using
$BlS48$-curve is $34778M$ \cite{NARDIEFO19}.

\section{Computational costs of the optimal Ate pairing using parallelizing elliptic net algorithm}\label{sectionparallel}
$\hspace*{1cm}$In this section, we parallelize the computation of
the optimal Ate pairing on the studied curves. The method we used is
stated in \cite[Section 4]{HiroOnu16}. This method helps to save one
multiplication in the addition step when a modified elliptic net is
considered. The method consists of defining a modified elliptic net
$W^{(1)}$ of $\widetilde{W}$ as
\begin{equation}
W^{(1)}(u,v)=\widetilde{W}(-1,1)^{uv}\widetilde{W}(u,v)~~ \forall u~
\forall v~\in \mathbb{Z}.
\end{equation}
$\hspace*{1cm}$For simplicity, we generalize the study by
considering a field
$\mathbb{F}_{p^{k}}$. We set $e=k/\delta$, where  $\delta$ denotes the degree of the twisted curve. We then have,\\
$\widetilde{W}(1,1)=1$,
$\widetilde{W}(-1,1)\in\mathbb{F}_{p^{k/2}}$,
$W^{(1)}(1,1)\in\mathbb{F}_{p^{k/2}}$, $W^{(1)}(2,-1)\in\mathbb{F}_{p^{k}}$, $W^{(1)}(-1,1)=1$.\\
A modified elliptic net is an elliptic net. One can see this
from Definition $\ref{def1}$.\\
The elliptic net $W^{(1)}$ satisfies the following relations:\\
{\scriptsize
    \begin{flushleft}
        $L_{1}:=~W^{(1)}(2k-3,0)=W^{(1)}(k,0)W^{(1)}(k-2,0)^{3}-W^{(1)}(k-3,0)W^{(1)}(k-1,0)^{3},$
    \end{flushleft}
    \begin{flushleft}
        $L_{2}:=~W^{(1)}(2k-2,0)=\frac{W^{(1)}(k-1,0)W^{(1)}(k+1,0)W^{(1)}(k-2,0)^{2}-W^{(1)}(k-1,0)W^{(1)}(k-3,0)W^{(1)}(k,0)^{2}}{W^{(1)}(2,0)},$
    \end{flushleft}
    \begin{flushleft}
        $L_{3}:=~W^{(1)}(2k-1,0)=W^{(1)}(k+1,0)W^{(1)}(k-1,0)^{3}-W^{(1)}(k-2,0)W^{(1)}(k,0)^{3},$
    \end{flushleft}
    \begin{flushleft}
        $L_{4}:=~W^{(1)}(2k,0)=\frac{W^{(1)}(k,0)W^{(1)}(k+2,0)W^{(1)}(k-1,0)^{2}-W^{(1)}(k,0)W^{(1)}(k-2,0)W^{(1)}(k+1,0)^{2}}{W^{(1)}(2,0)},$
    \end{flushleft}
    \begin{flushleft}
        $L_{5}:=~W^{(1)}(2k+1,0)=W^{(1)}(k+2,0)W^{(1)}(k,0)^{3}-W^{(1)}(k-1,0)W^{(1)}(k+1,0)^{3},$
    \end{flushleft}
    \begin{flushleft}
        $L_{6}:=~W^{(1)}(2k+2,0)=\frac{W^{(1)}(k+1,0)W^{(1)}(k+3,0)W^{(1)}(k,0)^{2}-W^{(1)}(k+1,0)W^{(1)}(k-1,0)W^{(1)}(k+2,0)^{2}}{W^{(1)}(2,0)},$
    \end{flushleft}
    \begin{flushleft}
        $L_{7}:=~W^{(1)}(2k+3,0)=W^{(1)}(k+3,0)W^{(1)}(k+1,0)^{3}-W^{(1)}(k,0)W^{(1)}(k+2,0)^{3},$
    \end{flushleft}
    \begin{flushleft}
        $L_{8}:=~W^{(1)}(2k+4,0)=\frac{W^{(1)}(k+2,0)W^{(1)}(k+4,0)W^{(1)}(k+1,0)^{2}-W^{(1)}(k+2,0)W^{(1)}(k,0)W^{(1)}(k+3,0)^{2}}{W^{(1)}(2,0)},$
    \end{flushleft}
    \begin{flushleft}
        $L_{9}:=~W^{(1)}(2k+5,0)=W^{(1)}(k+4,0)W^{(1)}(k+2,0)^{3}-W^{(1)}(k+1,0)W^{(1)}(k+3,0)^{3},$
    \end{flushleft}
    \begin{flushleft}
        $T_{1}:=~W^{(1)}(2k-1,1)=\frac{W^{(1)}(k+1,1)W^{(1)}(k-1,1)W^{(1)}(k-1,1)^{2}-W^{(1)}(k,0)W^{(1)}(k-2,0)W^{(1)}(k,1)^{2}}{W^{(1)}(1,1)},$
    \end{flushleft}
    \begin{flushleft}
        $T_{2}:=~W^{(1)}(2k,1)=W^{(1)}(k-1,1)W^{(1)}(k+1,1)W^{(1)}(k,0)^{2}-W^{(1)}(k-1,0)W^{(1)}(k+1,0)$\\$W^{(1)}(k,1)^{2},$
    \end{flushleft}
    \begin{flushleft}
        $T_{3}:=~W^{(1)}(2k+1,1)=W^{(1)}(k-1,1)W^{(1)}(k+1,1)W^{(1)}(k+1,0)^{2}-W^{(1)}(k,0)W^{(1)}(k+2,0)$\\$W^{(1)}(k,1)^{2},$
    \end{flushleft}
    \begin{flushleft}
        $T_{4}:=~W^{(1)}(2k+2,1)=\frac{W^{(1)}(k-1,1)W^{(1)}(k+1,1)W^{(1)}(k+2,0)^{2}-W^{(1)}(k+1,0)W^{(1)}(k+3,0)W^{(1)}(k,1)^{2}}{W^{(1)}(-2,1)}$.
    \end{flushleft}
}

The $\textbf{doubling step}$ consists of
calculating~$L_{1},~L_{2},~L_{2},~L_{4},~L_{5},~L_{6}$,
$L_{7},\\~L_{8},~T_{1},~ T_{2}$ and~$T_{3}$ whereas the addition
step consists of calculating
$L_{2},~L_{3},~L_{4},\\~L_{5},~L_{6},~L_{7},~L_{8},~L_{9},~T_{2},~
T_{3}$ and~$T_{4}$. The $\textbf{doubling step}$ in terms of
$W^{(1)}$ and the $\textbf{doubling step}$ in terms of
$\widetilde{W}$ have the same computational costs whereas the
$\textbf{addition step}$ in terms of $W^{(1)}$ entirely has 1
multiplication cost less than the computational cost of the
$\textbf{addition step}$ in terms of
$\widetilde{W}$. That is, with the new elliptic net, one multiplication is save.\\
$\hspace*{1cm}$In this work, we study the computational cost for
elliptic net algorithm with 4 and 8 processors.\\
Table 3 gives some notations and computational costs for elementary
factors and terms in the doubling and addition steps formulae,
in order to clearly present our parallel executions of optimal Ate pairings.\\
\begin{table}
    \centering {\scriptsize
        \begin{tabular}{|l|l|c|}
            \hline
            Operations & Values & Costs \\
            \hline
            $U_{1}:=W^{(1)}(k,0)W^{(1)}(k-2,0)$ & $U_{1}=W^{(1)}(k,0)W^{(1)}(k-2,0)$ & $M_{e}$ \\
            $U_{2}:=W^{(1)}(k-2,0)^{2}$ & $U_{2}=W^{(1)}(k-2,0)^{2}$& $S_{e}$ \\
            $U_{3}:=W^{(1)}(k-3,0)W^{(1)}(k-1,0)^{3}$ & $U_{3}=W^{(1)}(k-3,0)W^{(1)}(k-1,0)$ & $M_{e}$ \\
            $U_{4}:=W^{(1)}(k-1,0)^{2}$ & $U_{4}=W^{(1)}(k-1,0)^{2}$ & $S_{e}$ \\
            $U_{5}:=W^{(1)}(k-1,0)W^{(1)}(k+1,0)$&$U_{5}=W^{(1)}(k-1,0)W^{(1)}(k+1,0)$&$M_{e}$ \\
            $U_{6}:=W^{(1)}(k,0)^{2}$ & $U_{6}=W^{(1)}(k,0)^{2}$ & $S_{e}$ \\
            $U_{7}:=W^{(1)}(k,0)W^{(1)}(k+2,0)$&$U_{7}=W^{(1)}(k,0)W^{(1)}(k+2,0)$& $M_{e}$ \\
            $U_{8}:=W^{(1)}(k+1,0)^{2}$&$U_{8}=W^{(1)}(k+1,0)^{2}$&  $S_{e}$ \\
            $U_{9}:=W^{(1)}(k+1,0)W^{(1)}(k+3,0)$&$U_{9}=W^{(1)}(k+1,0)W^{(1)}(k+3,0)$&$M_{e}$ \\
            $U_{10}:=W^{(1)}(k+2,0)^{2}$&$U_{10}=W^{(1)}(k+2,0)^{2}$&$S_{e}$ \\
            $U_{11}:=W^{(1)}(k+2,0)W^{(1)}(k+4,0)$&$U_{11}=W^{(1)}(k+2,0)W^{(1)}(k+4,0)$&$M_{e}$\\
            $U_{12}:=W^{(1)}(k+3,0)^{2}$&$U_{12}=W^{(1)}(k+3,0)^{2}$&$S_{e}$\\
            $L_{1}:=U_{1}U_{2}-U_{3}U_{4}$ &$L_{1}=U_{1}U_{2}-U_{3}U_{4}$&$2M_{e}$\\
            $L_{2}:=\frac{1}{W(2,0)}(U_{5}U_{2}-U_{3}U_{6})$&$L_{2}=\frac{1}{W(2,0)}(U_{5}U_{2}-U_{3}U_{6})$& $2M_{e}$\\
            $L_{3}:=U_{5}U_{4}-U_{1}U_{6}$&$L_{3}=U_{5}U_{4}-U_{1}U_{6}$& $2M_{e}$ \\
            $L_{4}:=\frac{1}{W(2,0)}(U_{7}U_{4}-U_{1}U_{8})$&$L_{4}=\frac{1}{W(2,0)}(U_{7}U_{4}-U_{1}U_{8})$& $2M_{e}$\\
            $L_{5}:=U_{7}U_{6}-U_{5}U_{8}$&$L_{5}=U_{7}U_{6}-U_{5}U_{8}$&$2M_{e}$\\
            $L_{6}:=\frac{1}{W(2,0)}(U_{9}U_{6}-U_{5}U_{10})$&$L_{6}=\frac{1}{W(2,0)}(U_{9}U_{6}-U_{5}U_{10})$& $2M_{e}$\\
            $L_{7}:=U_{9}U_{8}-U_{7}U_{10}$&$L_{7}=U_{9}U_{8}-U_{7}U_{10}$& $2M_{e}$\\
            $L_{8}:=\frac{1}{W(2,0)}(U_{11}U_{8}-U_{7}U_{12})$&$L_{8}=\frac{1}{W(2,0)}(U_{11}U_{8}-U_{7}U_{12})$& $2M_{e}$\\
            $L_{9}:=U_{11}U_{10}-U_{9}U_{12}$&$L_{9}=U_{11}U_{10}-U_{9}U_{12}$& $2M_{e}$\\
            $V_{1}:=W^{(1)}(k+1,1)W^{(1)}(k-1,1)$&$V_{1}=W^{(1)}(k+1,1)W^{(1)}(k-1,1)$&$M_{k}$\\
            $V_{2}:=W^{(1)}(k,1)^{2}$&$V_{2}=W^{(1)}(k,1)^{2}$& $S_{k}$\\
            $X_{0}:=V_{1}U_{4}$& $X_{0}=V_{1}U_{4}$&$\delta M_{e}$\\
            $X_{1}:=V_{2}U_{1}$& $X_{1}=V_{2}U_{1}$&$\delta M_{e}$\\
            $X_{2}:=V_{1}U_{6}$& $X_{2}=V_{1}U_{6}$&$\delta M_{e}$\\
            $X_{3}:=V_{2}U_{5}$& $X_{3}=V_{2}U_{5}$&$\delta M_{e}$\\
            $X_{4}:=V_{1}U_{8}$& $X_{4}=V_{1}U_{8}$&$\delta M_{e}$\\
            $X_{5}:=V_{2}U_{7}$& $X_{5}=V_{2}U_{7}$&$\delta M_{e}$\\
            $X_{6}:=V_{2}U_{9}$& $X_{6}=V_{2}U_{9}$&$\delta M_{e}$\\
            $X_{7}:=V_{1}U_{10}$& $X_{7}=V_{1}U_{10}$&$\delta M_{e}$\\
            $Y_{1}:=X_{0}-X_{1}$& $Y_{1}=X_{0}-X_{1}$&$\cdots$\\
            $Y_{4}:=X_{6}-X_{7}$& $Y_{4}=X_{6}-X_{7}$&$\cdots$\\
            $T_{1}:=\frac{1}{W^{(1)}(1,1)}Y_{1}$&$T_{1}=\frac{1}{W^{(1)}(1,1)}Y_{1}$&$2M_{k/2}$\\
            $T_{2}:=X_{2}-X_{3}$& $T_{2}=X_{2}-X_{3}$&$\cdots$\\
            $T_{3}:=X_{4}-X_{5}$& $T_{3}=X_{4}-X_{5}$&$\cdots$\\
            $T_{4}:=\frac{1}{W^{(1)}(2,-1)}Y_{4}$& $T_{4}=\frac{1}{W^{(1)}(2,-1)}Y_{4}$&$M_{k}$\\
            \hline
        \end{tabular}\caption{Computational costs for some factors and terms in the doubling and addition steps formulae}.
    }

\end{table}

Table 4, Table 5, Table 6 and Table 7 give details for algorithm
with 4 and 8 processors respectively. For an embedding degree $k$
and a twist $\delta$ of a curve, we set
$e=k/\delta$.\\\\
\begin{table}
    \centering \scriptsize
    \begin{tabular}{|l|l|l|l|}
        \hline
        \multicolumn {2}{|c|}{processor 1} & \multicolumn {2}{c|}{Processor 2} \\
        \hline
        Operations & Costs &Operations & Costs \\
        $U_{1},~U_{3},~U_{5}$& $3(M_{e})$ &$U_{1},~U_{5},~U_{6}$& $\cdots$ \\
        $U_{2},~U_{4},~U_{6}$& $3(S_{e})$ &$U_{7}//~U_{8}$& $M_{e}//~S_{e}$ \\
        $V_{1}$& $M_{k}$ &$V_{2}$& $S_{k}$ \\
        $L_{1}:=U_{1}U_{2}-U_{3}U_{4}$ & $2M_{e}$ & $L_{3}:=U_{5}U_{4}-U_{1}U_{6}$ & $2M_{e}$ \\
        $L_{2}:=\frac{1}{W^{(1)}(2,0)}(U_{5}U_{2}-U_{3}U_{6})$ &$2M_{e}$ &$L_{4}:=\frac{1}{W^{(1)}(2,0)}(U_{7}U_{4}-U_{1}U_{8})$ &$2M_{e}$ \\
        $X_{0}:=V_{1}U_{4}$  & $\delta M_{e}$ &$X_{1}:=V_{2}U_{1}$  & $\delta M_{e}$  \\
        $X_{2}:=V_{1}U_{6}$  &  $\delta M_{e}$  &$X_{3}:=V_{2}U_{5}$  &  $\delta M_{e}$  \\
        \hline
        &&&\\
        & $(7+2\delta)M_{e}+3S_{e}+M_{k}$ &   & $(5+2\delta)M_{e}+S_{e}+S_{k}$ \\
        \hline

        \multicolumn {2}{|c|}{processor 3} & \multicolumn {2}{c|}{Processor 4} \\
        \hline
        Operations & Costs &Operations & Costs \\
        $U_{5},~U_{6},~U_{7},~U_{8}$& $\cdots$ &$U_{7},~U_{8},~U_{9},~U_{10}$& $\cdots$ \\
        $U_{9}//~U_{10}$& $M_{e}//~S_{e}$ &$U_{11}//~U_{12}$& $M_{e}//~S_{e}$ \\
        $V_{1}$& $\cdots$ & $V_{2}$& $\cdots$  \\
        $L_{5}:=U_{7}U_{6}-U_{5}U_{8}$ & $2M_{e}$ &$L_{7}:=U_{9}U_{8}-U_{7}U_{10}$ &  $2M_{e}$\\
        $L_{6}:=\frac{1}{W^{(1)}(2,0)}(U_{9}U_{6}-U_{5}U_{10})$ &$2M_{e}$ &$L_{8}:=\frac{1}{W^{(1)}(2,0)}(U_{11}U_{8}-U_{7}U_{12})$ &$2M_{e}$\\
        $X_{4}:=V_{1}U_{8}$  & $\delta M_{e}$ & $X_{5}:=V_{2}U_{7}$ &  $\delta M_{e}$\\
        $X_{0},~X_{1},~X_{2},~X_{4}$  &  $\cdots$  &$X_{4},~Y_{1},~T_{3}$ & $\cdots$ \\
        $Y_{1},~T_{2}$& $\cdots$ &$T_{1}:=\frac{1}{W^{(1)}(1,1)}Y_{1}$ & $2 M_{k/2}$ \\

        \hline
        &$(5+\delta)M_{e}+S_{e}$    &     & $(5+\delta)M_{e}+S_{e}+2M_{k/2}$  \\
        \hline

        \multicolumn {4}{|c|}{Computational cost of the longest path~~$(7+2\delta)M_{e}+3S_{e}+M_{k}$} \\
        \hline
    \end{tabular}\caption{Computational costs of the Doubling step for algorithm with 4 processors }.
    $$
    $$

    \begin{tabular}{|l|l|l|l|}
        \hline
        \multicolumn {2}{|c|}{processor 1} & \multicolumn {2}{c|}{Processor 2} \\
        \hline
        Operations & Costs &Operations & Costs \\
        $U_{1},~U_{3},~U_{5}$& $3(M_{e})$ &$U_{1},~U_{4},~U_{5},~U_{6},~U_{8}$& $\cdots$ \\
        $U_{2},~U_{4},~U_{6},~U_{8}$& $4(S_{e})$ &$U_{7}$& $M_{e}$ \\
        $V_{1}$& $M_{k}$ &$V_{2}$& $S_{k}$ \\
        $L_{2}:=\frac{1}{W^{(1)}(2,0)}(U_{5}U_{2}-U_{3}U_{6})$ & $2M_{e}$ & $L_{4}:=\frac{1}{W^{(1)}(2,0)}(U_{7}U_{4}-U_{1}U_{8})$ & $2M_{e}$ \\
        $L_{3}:=U_{5}U_{4}-U_{1}U_{6}$   &$2M_{e}$ &$L_{5}:=U_{7}U_{6}-U_{5}U_{8}$ &$2M_{e}$ \\
        $X_{2}:=V_{1}U_{6}$  & $\delta M_{e}$ &$X_{3}:=V_{2}U_{5}$  & $\delta M_{e}$  \\
        $X_{4}:=V_{1}U_{8}$  &  $\delta M_{e}$  &$X_{5}:=V_{2}U_{7}$  &  $\delta M_{e}$  \\
        &   &$X_{2},~X_{4}$  &  $\cdots$  \\
        &   &$T_{2},~T_{3}$  &  $\cdots$  \\
        \hline
        &&&\\
        & $(7+2\delta)M_{e}+4S_{e}+M_{k}$ &   & $(5+2\delta)M_{e}+S_{k}$ \\
        \hline

        \multicolumn {2}{|c|}{processor 3} & \multicolumn {2}{c|}{Processor 4} \\

        \hline
        Operations & Costs &Operations & Costs \\
        $U_{5},~U_{6},~U_{7},~U_{8}$& $\cdots$ &$U_{7},~U_{8},~U_{9},~U_{10}$& $\cdots$ \\
        $U_{9}//~U_{10}$& $M_{e}//~S_{e}$ &$U_{11}//~U_{12}$& $M_{e}//~S_{e}$ \\
        $V_{1}$& $\cdots$ & $V_{1}$& $\cdots$  \\
        $L_{6}:=\frac{1}{W^{(1)}(2,0)}(U_{9}U_{6}-U_{5}U_{10})$ & $2M_{e}$ &$L_{8}:=\frac{1}{W^{(1)}(2,0)}(U_{11}U_{8}-U_{7}U_{12})$ &  $2M_{e}$\\
        $L_{7}:=U_{9}U_{8}-U_{7}U_{10}$ &$2M_{e}$ &$L_{9}:=U_{11}U_{10}-U_{9}U_{12}$ &$2M_{e}$\\
        $X_{6}:=V_{2}U_{9}$  & $\delta M_{e}$ & $X_{7}:=V_{1}U_{10}$ &  $\delta M_{e}$\\
        $X_{7},~Y_{4}$  &  $\cdots$  &     & \\
        $T_{4}:=\frac{1}{W^{(1)}(-2,1)}Y_{4}$& $2 M_{k/2}$ &    & \\

        \hline
        &$(5+\delta)M_{e}+S_{e}+2M_{k/2}$    &     & $(5+\delta)M_{e}+S_{e}$  \\
        \hline

        \multicolumn {4}{|c|}{Computational cost of the longest path~~$(7+2\delta)M_{e}+4S_{e}+M_{k}$} \\
        \hline
    \end{tabular}\caption{Computational costs of the Addition step for algorithm with 4 processors }.

\end{table}

\begin{table}
    \centering \scriptsize
    \begin{tabular}{|l|l|l|l|}
        \hline
        \multicolumn {2}{|c|}{processor 1} & \multicolumn {2}{c|}{Processor 2} \\
        \hline
        Operations & Costs &Operations & Costs \\
        $U_{1},~U_{3}$& $2(M_{e})$ &$U_{2},~U_{3}$& $\cdots$ \\
        $U_{2},~U_{4}$& $2(S_{e})$ &$U_{5}//~U_{6}$& $M_{e}//~S_{e}$ \\
        $V_{1}$& $M_{k}$ &$V_{1}$& $\cdots$ \\
        $L_{1}:=U_{1}U_{2}-U_{3}U_{4}$ & $2M_{e}$ & $L_{2}:=\frac{1}{W^{(1)}(2,0)}(U_{5}U_{2}-U_{3}U_{6})$ & $2M_{e}$ \\
        $X_{0}:=V_{1}U_{4}$  & $\delta M_{e}$ &$X_{2}:=V_{1}U_{6}$  & $\delta M_{e}$  \\
        &  & &   \\
        \hline
        &&&\\
        & $(4+\delta)M_{e}+2S_{e}+M_{k}$ &   & $(3+\delta)M_{e}+S_{e}$ \\
        \hline

        \multicolumn {2}{|c|}{processor 3} & \multicolumn {2}{c|}{Processor 4} \\
        \hline
        Operations & Costs &Operations & Costs \\
        $U_{1},~U_{4},~U_{5},~U_{6}$& $\cdots$ &$U_{1},~U_{4}$& $\cdots$ \\
        $V_{2}$    &$S_{k}$    &$U_{7}//~U_{8}$& $M_{e}//~S_{e}$ \\
        $L_{3}:=U_{5}U_{4}-U_{1}U_{6}$   &$2M_{e}$ &$L_{4}:=\frac{1}{W^{(1)}(2,0)}(U_{7}U_{4}-U_{1}U_{8})$ & $2M_{e}$ \\
        $X_{1}:=V_{2}U_{1}$  & $\delta M_{e}$ & $X_{0},~X_{1}$  &  $\cdots$   \\
        $X_{4},~X_{5}$  &  $\cdots$  &$Y_{1}$& $\cdots$ \\
        $Y_{3}$& $\cdots$     & $T_{1}:=\frac{1}{W^{(1)}(1,1)}Y_{1}$ & $2 M_{k/2}$  \\

        \hline
        &$(2+\delta)M_{e}+S_{k}$    &     & $3M_{e}+S_{e}+2M_{k/2}$  \\
        \hline

        \multicolumn {2}{|c|}{processor 5} & \multicolumn {2}{c|}{Processor 6} \\
        \hline
        Operations &Costs& Operations& Costs\\
        $U_{5},~U_{6},~U_{7},~U_{8}$& $\cdots$ &$U_{5},~U_{6}$& $\cdots$ \\
        $L_{5}:=U_{7}U_{6}-U_{5}U_{8}$ & $2M_{e}$  &$U_{9}//~U_{10}$& $M_{e}//~S_{e}$ \\
        $V_{2}$& $\cdots$ &$L_{6}:=\frac{1}{W^{(1)}(2,0)}(U_{9}U_{6}-U_{5}U_{10})$ & $2M_{e}$ \\
        $X_{3}:=V_{2}U_{5}$  & $\delta M_{e}$  &  &  \\
        &   &   &\\
        \hline
        & $(2+\delta)M_{e}$ &   & $3M_{e}+S_{e}$ \\
        \hline

        \multicolumn {2}{|c|}{processor 7} & \multicolumn {2}{c|}{Processor 8} \\
        \hline
        Operations & Costs &Operations & Costs \\
        $U_{7},~U_{8},~U_{9},~U_{10}$& $\cdots$ &$U_{7},~U_{8}$& $\cdots$ \\
        $L_{7}:=U_{9}U_{8}-U_{7}U_{10}$&$2M_{e}$ &$U_{11}//~U_{12}$& $M_{e}//~S_{e}$ \\
        $V_{1}$& $\cdots$ & $L_{8}:=\frac{1}{W^{(1)}(2,0)}(U_{11}U_{8}-U_{7}U_{12})$ & $2M_{e}$ \\
        $X_{4}:=V_{1}U_{8}$  & $\delta M_{e}$&$V_{2}$& $\cdots$\\
        &    & $X_{2},~X_{3}$  &  $\cdots$ \\
        &   & $T_{2}$ &  $\cdots$\\
        \hline
        &$(2+\delta)M_{e}$    &     & $3M_{e}+S_{e}$  \\
        \hline

        \multicolumn {4}{|c|}{Computational cost of the longest path~~$(4+\delta)M_{e}+2S_{e}+M_{k}$} \\
        \hline
    \end{tabular}\caption{Computational costs of the Doubling step for algorithm with 8 processors }.

\end{table}

\begin{table}
    \centering \scriptsize
    \begin{tabular}{|l|l|l|l|}
        \hline
        \multicolumn {2}{|c|}{processor 1} & \multicolumn {2}{c|}{Processor 2} \\
        \hline
        Operations & Costs &Operations & Costs \\
        $U_{3},~U_{5}$& $2(M_{e})$ &$U_{5},~U_{6}$& $\cdots$ \\
        $U_{2},~U_{6}$& $2(S_{e})$ &$U_{1}//~U_{4}$& $M_{e}//~S_{e}$ \\
        $V_{1}$& $M_{k}$ &$V_{2}$& $S_{k}$ \\
        $L_{2}:=\frac{1}{W^{(1)}(2,0)}U_{5}U_{2}-U_{3}U_{6}$ & $2M_{e}$ & $L_{3}:=U_{5}U_{4}-U_{1}U_{6}$ & $2M_{e}$ \\
        $X_{2}:=V_{1}U_{6}$  & $\delta M_{e}$ &$X_{3}:=V_{2}U_{5}$  & $\delta M_{e}$  \\
        &  & &   \\
        \hline
        &&&\\
        & $(4+\delta)M_{e}+2S_{e}+M_{k}$ &   & $(3+\delta)M_{e}+S_{e}+S_{k}$ \\
        \hline

        \multicolumn {2}{|c|}{processor 3} & \multicolumn {2}{c|}{Processor 4} \\
        \hline
        Operations & Costs &Operations & Costs \\
        $U_{1},~U_{4}$& $\cdots$ &$U_{5},~U_{6},~U_{7},~U_{8}$& $\cdots$ \\
        $U_{7}//~U_{8}$& $M_{e}//~S_{e}$ & $L_{5}:=U_{7}U_{6}-U_{5}U_{8}$ & $2M_{e}$ \\
        $L_{4}:=\frac{1}{W^{(1)}(2,0)}(U_{7}U_{4}-U_{1}U_{8})$ & $2M_{e}$ & $V_{2}$& $\cdots$   \\
        $V_{1}$& $\cdots$ & $X_{5}:=V_{2}U_{7}$  & $\delta M_{e}$ \\
        $X_{4}:=V_{1}U_{8}$  & $\delta M_{e}$ &   & \\
        \hline
        &$(3+\delta)M_{e}+S_{e}$    &     & $(2+\delta)M_{e}$  \\
        \hline

        \multicolumn {2}{|c|}{processor 5} & \multicolumn {2}{c|}{Processor 6} \\
        \hline
        Operations &Costs& Operations& Costs\\
        $U_{5},~U_{6}$& $\cdots$ &$U_{7},~U_{8},~U_{9},~U_{10}$& $\cdots$ \\
        $U_{9}//~U_{10}$& $M_{e}//~S_{e}$ &$L_{7}:=U_{9}U_{8}-U_{7}U_{10}$ & $2M_{e}$  \\
        $L_{6}:=\frac{1}{W^{(1)}(2,0)}(U_{9}U_{6}-U_{5}U_{10})$ & $2M_{e}$  &$V_{1}$& $\cdots$ \\
        $V_{2}$& $\cdots$ &$X_{7}:=V_{1}U_{10}$  & $\delta M_{e}$  \\
        $X_{6}:=V_{2}U_{9}$  & $\delta M_{e}$  &  &  \\
        &   &   &\\
        \hline
        & $(3+\delta)M_{e}+S_{e}$ &   & $(2+\delta)M_{e}$ \\
        \hline

        \multicolumn {2}{|c|}{processor 7} & \multicolumn {2}{c|}{Processor 8} \\
        \hline
        Operations & Costs &Operations & Costs \\
        $U_{7},~U_{8}$& $\cdots$ &$U_{9},~U_{10},~U_{11},~U_{12}$& $\cdots$ \\
        $U_{11}//~U_{12}$& $M_{e}//~S_{e}$ &$L_{9}:=U_{11}U_{10}-U_{9}U_{12}$&$2M_{e}$\\
        $L_{8}:=\frac{1}{W^{(1)}(2,0)}(U_{11}U_{8}-U_{7}U_{12})$&$2M_{e}$ &$Y_{4}$& $\cdots$ \\
        $X_{2},~X_{3},~X_{4},~X_{5},~X_{6},~X_{7}$  &  $\cdots$  & $T_{4}:=\frac{1}{W^{(1)}(2,-1)}Y_{4}$ & $M_{k}$ \\
        $T_{2},~T_{3},~Y_{4}$  &  $\cdots$ &$V_{2}$& $\cdots$\\
        \hline
        &$3M_{e}+S_{e}$    &     & $2M_{e}+M_{k}$  \\
        \hline

        \multicolumn {4}{|c|}{Computational cost of the longest path~~$(4+\delta)M_{e}+2S_{e}+M_{k}$} \\
        \hline
    \end{tabular}\caption{Computational costs of the Addition step for algorithm with 8 processors }.

\end{table}

\subsection{Notation and cost of the arithmetic in finite field}
In this work, $M, S$ and $I$ denote the cost of one multiplication,
one squaring, one inversion in the finite field $\mathbb{F}_{p}$
respectively. $M_{i}, S_{i}$ and $I_{i}$ denote the computation
costs for one multiplication, one squaring and one inversion in the
finite extension field $\mathbb{F}_{p^{i}}$ of $\mathbb{F}_{p}$.\\
From \cite{Aranha11}, one can have the following costs:\\
$M_{2}=3M$, $S_{2}=\frac{2}{3}M_{2}$ $M_{3}=6M$, $S_{3}=5M$,
$M_{4}=9M$, $S_{4}=6M$, $M_{6}=18M$, $M_{8}=27M$, $S_{8}=18M$,
$M_{9}=36M$, $S_{9}=25M$, $M_{16}=81M$, $S_{16}=54M$, $M_{18}=108M$,
$S_{18}=55M$, $M_{24}=162M$, $S_{24}=108M$, $M_{48}=486M$,
$S_{48}=324M$,. $I_{6}=37M+I$.  1 $p$ and $p^2-$ Frobenius in
$\mathbb{F}_{p^{12}}$ are respectively $10M$ and $15M$, 1 $p$ and
$p^3-$ Frobenius in $\mathbb{F}_{p^{16}}$ are $15M$ each.\\

\subsection{Computational costs of optimal Ate pairing on the studied curves.}
In this subsection, we provide the computational costs of the part
without the final exponentiation of the optimal Ate pairing at 128,
192 and
256-bit security level on our chosen curves.\\
\subsubsection{Computational costs at 128-bit security level.}

\begin{itemize}
    \item Case of BN-curve. Let
    $f_{1}=\widetilde{W}(6x+2,1)\cdot\mathbb{L}_{1}\cdot\mathbb{L}_{2}$
    be the part without the final exponentiation of the optimal Ate
    pairing on a BN-curve. Then the value $x=2^{114}+2^{101}-2^{14}-1$
    from Table 2 gives
    $6x+2=2^{116}+2^{115}+2^{103}+2^{102}-2^{16}-2^{15}-2^{2}$.
    $\widetilde{W}(6x+2,1)$ costs 116 Doubling steps and 6 addition
    steps.\\
    \begin{enumerate}
        \item For 4 processors,the doubling step costs
        $(7+2\times6)M_{2}+3S_{2}+M_{12}$ and the addition step costs
        $(7+2\times6)M_{2}+4S_{2}+M_{12}$, that is $117M$ and $119M$
        respectively, contrary to $108M$ and $126M$ from
        \cite{HiroOnu16}.\\
        Then, $\widetilde{W}(6x+2,1)$ costs 116 Doubling steps and 6
        addition steps, that is $116(117M)+6(119M)=14286M$. $[6x+2]Q$,
        $[p]Q$ and $\mathbb{L}_{1}$ together cost $2M_{2}+2I_{2}$, 2
        $p$-Frobenius in $\mathbb{F}_{p^{12}}$ and $76M$, that is $100M+I$.
        $[6x+2+p]Q$, $[-p^{2}]Q$ and $\mathbb{L}_{2}$ cost $340M+2I$, 2
        $p^{2}$-Frobenius in $\mathbb{F}_{p^{12}}$ and $222M$, that is
        $577M+2I$. $f_{1}$ costs
        $14286M+(100M+I)+(577M+2I)+2M_{12}$, that $15071M+3I$.\\
        \item For 8 processors, the doubling step and addition step cost
        $(4+6)M_{2}+2S_{2}+M_{12}$ and $(4+6)M_{2}+2S_{2}+M_{12}$
        respectively. That is $88M$ and $88M$ respectively, contrary to
        $90M$ and $90M$ from \cite{HiroOnu16}. $\widetilde{W}(6x+2,1)$
        costs $116(88M)+6(88M)=10736M$. $f_{1}$ then costs $11521M+3I$.\\
    \end{enumerate}
    \item Case of BLS12-curve. Let $f_{2}=\widetilde{W}(x,1)$ be the part
    without the final exponentiation of the optimal Ate pairing on a
    BLS12-curve. For the selected parameter $x=-2^{77}+2^{50}+2^{33}$,
    $\widetilde{W}(x,1)$ costs 77 Doubling steps and 2 addition steps.\\
    \begin{enumerate}
        \item For 4 processors, the doubling step and the addition step cost
        $(7+2\times6)M_{2}+3S_{2}+M_{12}$ and
        $(7+2\times6)M_{2}+4S_{2}+M_{12}$ respectively, that is $117M$ and
        $119M$. $\widetilde{W}(x,1)$ then costs  $77(117M)+2(119M)=9247M$.\\
        \item For 8 processors, the doubling step and addition step cost
        $(4+6)M_{2}+2S_{2}+M_{12}$ and $(4+6)M_{2}+2S_{2}+M_{12}$. That is
        $88M$ and $88M$ respectively, contrary to $90M$ and $90M$ from
        \cite{HiroOnu16}. $\widetilde{W}(x,1)$ costs
        $77(88M)+2(88M)=6952M$. $f_{2}$ then costs $6952M$.\\
    \end{enumerate}
    \item Case of KSS16-curve. Let
    $f_{3}=(\widetilde{W}(x,1)\cdot\ell_{1})^{p^{3}}\cdot\ell_{2}$ be
    the part without the final exponentiation of the optimal Ate pairing
    on a KSS16-curve. For the selected parameter
    $x=2^{35}-2^{32}-2^{18}+2^{8}+1$,
    $\widetilde{W}(x,1)$ costs 35 Doubling steps and 4 addition steps.\\
    \begin{enumerate}
        \item For 4 processors, the doubling step and the addition step cost
        $(7+2\times4)M_{4}+3S_{4}+M_{16}$ and
        $(7+2\times4)M_{4}+4S_{4}+M_{16}$ respectively, that is $234M$ and
        $240M$. $\widetilde{W}(x,1)$ costs $35(234M)+4(240M)=9150M$. $[x]Q$,
        $[p]Q$, $\ell_{1}$ and $(\widetilde{W}(x,1)\cdot\ell_{1})^{p^{3}}$
        together cost $395M+2I$. $\ell_{2}$ costs $92M$ and $f_{3}$ then
        costs $9637M+2I$.\\
        \item For 8 processors, the doubling step and addition step cost
        $(4+4)M_{4}+2S_{4}+M_{16}$ and $(4+4)M_{4}+2S_{4}+M_{16}$. That is
        $165M$ and $165M$ respectively.  $\widetilde{W}(x,1)$ costs
        $35(165M)+4(165M)=6435M$. $f_{3}$ then costs $6922M+2I$.\\
    \end{enumerate}

\end{itemize}

\subsubsection{Computational costs at 192-bit security level.}

\begin{itemize}

    \item Case of BLS24-curve. Let $f_{5}=\widetilde{W}(x,1)$ be the part
    without the final exponentiation of the optimal Ate pairing on a
    BLS24-curve. For the selected parameter
    $x=-2^{56}-2^{43}+2^{9}-2^{6}$,
    $\widetilde{W}(x,1)$ costs 56 Doubling steps and 3 addition steps.\\
    \begin{enumerate}
        \item For 4 processors, the doubling step and the addition step cost
        $(7+2\times6)M_{4}+3S_{4}+M_{24}$ and
        $(7+2\times6)M_{4}+4S_{4}+M_{24}$ respectively. That is $351M$ and
        $357M$ respectively. $f_{5}=\widetilde{W}(x,1)$ costs
        $56(351M)+3(357M)=20727M$.
        \item For 8 processors, the doubling step and addition step cost $(4+6)M_{4}+2S_{4}+M_{24}$ and
        $(4+6)M_{4}+2S_{4}+M_{24}$ respectively $264M$ and $264M$
        respectively. $f_{5}=\widetilde{W}(x,1)$ costs
        $56(264M)+3(264M)=15576M$
    \end{enumerate}

\end{itemize}

\subsubsection{Computational costs at 256-bit security level.}

\begin{itemize}
    \item Case of BLS24-curve. Let $f_{5}=\widetilde{W}(x,1)$ be the part
    without the final exponentiation of the optimal Ate pairing on a
    BLS24-curve. For the selected parameter
    $x=-2^{103}-2^{101}+2^{68}+2^{50}$,
    $\widetilde{W}(x,1)$ costs 103 Doubling steps and 3 addition steps.\\
    \begin{enumerate}
        \item For 4 processors, the doubling step and the addition step cost
        $351M$ and $357M$ respectively. $f_{5}=\widetilde{W}(x,1)$ costs
        $103(351M)+3(357M)=37224M$.
        \item For 8 processors, the doubling step and addition step cost $264M$ and $264M$
        respectively. $f_{5}=\widetilde{W}(x,1)$ costs
        $103(264M)+3(264M)=27984M$.\\
    \end{enumerate}

    \item  Case of BLS48-curve. Let $f_{6}=\widetilde{W}(x,1)$ be the part
    without the final exponentiation of the optimal Ate pairing on a
    BLS48-curve. For the selected parameter
    $x=2^{32}-2^{18}-2^{10}-2^{4}$,
    $\widetilde{W}(x,1)$ costs 32 Doubling steps and 3 addition steps.\\
    \begin{enumerate}
        \item For 4 processors, the doubling step and the addition step cost
        $(7+2\times6)M_{8}+3S_{8}+M_{48}$ and
        $(7+2\times6)M_{8}+4S_{8}+M_{48}$ respectively. That is $1053M$ and
        $1071M$ respectively. $f_{6}=\widetilde{W}(x,1)$ costs
        $32(1053M)+3(1071M)=36909M$.
        \item For 8 processors, the doubling step and addition step cost $(4+6)M_{8}+2S_{8}+M_{48}$ and
        $(4+6)M_{4}+2S_{4}+M_{24}$ respectively $264M$ and $264M$
        respectively. $f_{6}=\widetilde{W}(x,1)$ costs
        $32(792M)+3(792M)=27720M$.\\
    \end{enumerate}
%

\end{itemize}

\begin{table}
    \centering \scriptsize
    \begin{tabular}{|c|l|l|l|l|}
        \hline
        \multicolumn {5}{|c|}{Level 128} \\
        \hline
        ///////////// &Type of comp.& $f_{1}$ (BN-curve)&$f_{2}$ (BLS12-curve)&$f_{3}$ (KSS16-curve)  \\
        \hline
        Number of Proc.  &  Miller loop \cite{RazDuq18} & $12068M$ & $7708M$ & $7534M$ \\
        \hline
        4  &Elliptic net& $15071M+3I$ & $9247M$ & $9637M+2I$ \\
        \hline
        8&Elliptic net& $11521M+3I$  & $6952M$  & $6922M+2I$ \\
        \hline
    \end{tabular}\caption{Computational costs of the path without the final exponentiation of the optimal Ate pairing
        for the BN, BLS12 and KSS16 curves for 4 and 8 processors }.
    $$
    $$
    \begin{tabular}{|c|l|l|}
        \hline
        \multicolumn {3}{|c|}{Level 192} \\
        \hline
        ///////////// &Type of comp. & $f_{5}$ (BLS24-curve) \\
        \hline
        Number of Proc.  &  Miller loop & $19474M$  \\
        \hline
        4  &Elliptic net&  $20727M$  \\
        \hline
        8&Elliptic net&  $15576M$  \\
        \hline
    \end{tabular}\caption{Computational costs of the path without the final exponentiation of the optimal Ate pairing
        for BLS24 curves with 4 and 8 processors. }
    $$
    $$
    \begin{tabular}{|c|l|l|l|}
        \hline
        \multicolumn {4}{|c|}{Level 256} \\
        \hline
        ///////////// &Type of comp.&$f_{5}$ (BLS24-curve)&$f_{6}$ (BLS48-curve) \\
        \hline
        Number of Proc.  &  Miller loop& $35360M$ & $34778M$ \\
        \hline
        4  &Elliptic net& $37224M$ & $36909M$  \\
        \hline
        8&Elliptic net&  $27984M$  & $27720M$ \\
        \hline
    \end{tabular}\caption{Computational costs of the path without the final exponentiation of the optimal Ate pairing
        for the BLS24 and BLS48 curves for 4 and 8 processors. }
\end{table}
\newpage
\subsection{Comparison} $\hspace*{1cm}$One can see from Table 8, 9 and 10
that our parallel execution of the Elliptic Net Algorithm with 4
processors is less faster than the Miller method. But, for parallel
execution with 8 processors, the elliptic net method becomes more
faster.  The parallel execution of BLS48-curve with 8 processors is
more faster that the Miller method and when observing the results in
the last line of Table 10, the parallel execution of the optimal Ate
pairing on BLS48-curves is more faster than the one on BLS24 curves.

\section{Conclusion}\label{conclusion} In this work, we have computed the optimal Ate
pairing on BN, BLS12, KSS16, BLS24 and  BLS48 curves in terms of
elliptic nets associated twisted corresponding curves. We have given
the computational costs of the path without the final exponentiation
of the optimal Ate pairings  with 4 and 8 processors. We have seen
that the parallel execution with 8 processors give faster results
with elliptic nets, compared to Miller method.

\section*{Acknowledgment}
After the review
\bibliographystyle{unsrt}
\bibliography{bibliographynar6}
\end{document}